\documentclass{article}
\usepackage{caption}
\usepackage{amssymb}
\usepackage{amsmath}
\usepackage{amsfonts}
\usepackage{graphicx}
\usepackage{epstopdf}
\usepackage{placeins}
\usepackage{color}

\newcommand{\beq}{\begin{eqnarray}}
\newcommand{\eeq}{\end{eqnarray}}

\begin{document}

\begin{center}
\vspace{1cm}{\Large {\bf On Universal  Quantum Dimensions}}

\vspace{1cm}{\bf R.L.Mkrtchyan} \footnote{ E-mail: mrl55@list.ru} 

\vspace{1cm}

{\it Yerevan Physics Institute,}

{\it 2 Br. Alikhanian  Str., Yerevan 0036, Armenia}
\end{center}

\vspace{1cm}

\begin{abstract}
We represent in the universal form restricted one-instanton partition function of supersymmetric Yang-Mills theory. 
It is based on the derivation of universal expressions for quantum dimensions (universal characters) of Cartan powers of adjoint and  some other series of irreps of simple Lie algebras. These formulae also provide a proof of formulae for universal quantum dimensions for low-dimensional representations, needed in derivation of universal knot polynomials (i.e. colored Wilson averages of Chern-Simons theory on 3d sphere). As a check of the (complicated) formulae for universal quantum dimensions we prove numerically Deligne's hypothesis on universal characters for symmetric cube of adjoint representation.

Keywords: supersymmetric gauge theories, quantum dimensions, universal Lie algebra, instanton sums,  arXiv:1610.09910.

MSC classes: 17B20, 17B37

\end{abstract}


\section{Introduction}

This paper is the next step in realization of program of representation of partition functions and observables in gauge theories (and other simple-Lie-algebras-based theories) in universal form.
 
Universal formulae may appear in any theory, based on the simple Lie algebras. Most important example  is the Chern-Simons theory - see e.g. \cite{MV,M13} for universal expression for partition function and \cite{MMM,MM} for universal knot polynomials (Wilson loops). Universal formulae appear to be relevant also for refined Chern-Simons theory \cite{KS}, non-perturbative gauge/string duality \cite{KM} and Diophantine classification of simple Lie algebras \cite{M1,KhM1}.  Universality approach reveals \cite{KhM} details of behavior of $SU(N)$ Chern-Simons partition function under $N \rightarrow -N$ duality, and, more generally, under permutations of universal parameters (Vogel's symmetry). 

The main aim of present paper is derivation of universal formula for one-instanton contribution to partition function of Yang-Mills theory with extended supersymmetry (see Section \ref{s:1in}). This will be the first appearance of universal formula in four-dimensional Yang-Mills theory. Another aim is derivation  of formulae for quantum dimensions, used in \cite{MMM} for calculation of universal expressions for some knot polynomials (Wilson averages in Chern-Simons theory). 

Both results are based on the derivation (see Sections \ref{s:2}, \ref{s:3}) of new formulae for universal quantum dimensions (sometimes called also universal characters) of some series of irreps, appearing in decomposition of powers of adjoint representation. These formulae themselves can be considered as universal form of some observables - namely, unknot Wilson loops in corresponding representations.

The notion of universality in simple Lie algebras  was  introduced by Vogel in the paper "Universal Lie Algebra" \cite{V},  particularly aimed on creation of the universal model for all simple Lie algebras. In this approach different algebras would appear by specializing (homomorphic mapping) ring of coefficients of universal Lie algebra into that appropriate to a given simple Lie algebra.  Ring of coefficients of the universal Lie algebra is the ring of one-variable polynomials over the ring $\Lambda$ of antisymmetric three-leg Jacobi diagrams, introduced by Vogel in \cite{V}. It appears that one indeed can obtain in this way all simple Lie algebras (and superalgebras). Moreover, as shown in \cite{R},  some other algebras, namely Kapranov's algebra for $K_3$ manifold, can be obtained in this way, with appropriate choice of ring of coefficients. 

In this way one can try to develop the representation theory of universal Lie algebra, which particularly will contain the representation theory of powers of adjoint representation of simple Lie algebras. This was done in \cite{V}   for square and cube of adjoint representation. They are decomposed in universal form into simple modules, under some assumptions. When specialized to any simple Lie algebra, this decomposition gives true decomposition for that algebra. Formulae of \cite{V} for dimensions, as well as eigenvalues of second and fourth Casimir operators on that representations, are first examples of universal expressions. 

 For example, universal formula for dimension of simple Lie algebra $\mathfrak {g}$ is: 

\begin{eqnarray}\label{dim}
dim (\mathfrak {g}) &=& -\frac{(\alpha+2\beta+2\gamma)(\beta+2\alpha+2\gamma)(\gamma+2\alpha+2\beta)}{\alpha\beta\gamma}
\end{eqnarray}

Here $\alpha,\beta, \gamma$ are solutions of the cubic equation in ring $\Lambda(x)$, i.e. they belong to  some cubic extension of $\Lambda(x)$:

\begin{eqnarray}
\psi^3-t\psi^2+s\psi-p=0
\end{eqnarray}
where $t, s, p$ are certain elements of the ring, i.e. some combinations of Jacobi diagrams. According to Vieta's formulae one has

\begin{eqnarray}\label{}
t&=&\alpha+\beta+\gamma \\
s&=& \alpha\beta+\beta \gamma +\alpha\gamma \\
p&=&\alpha\beta\gamma
\end{eqnarray}
where  cubic equation's roots $\alpha,\beta,\gamma$ are called universal, or Vogel's parameters. They are relevant up to the rescaling and permutations \cite{V}, so belong to two-dimensional projective plane, factorized over permutations of homogeneous parameters $\alpha,\beta,\gamma$, which is called Vogel's plane.

When specialized for complex simple Lie algebras, these parameters get values from Vogel's table \ref{tab:V1}. The same table in other form (\ref{tab:V2}) reveals that not only  all orthogonal algebras belong to one line in Vogel's plane, but that the same line contains all symplectic algebras. Another similar and unexpected observation is that all exceptional algebras belong to one line, which we call exceptional, or Deligne line \cite{Del,DM}. 

\begin{table}[h]   
\caption{Vogel's parameters for simple Lie algebras}     
\begin{tabular}{|c|c|c|c|c|c|}
\hline
Root system & Lie algebra  & $\alpha$ & $\beta$ & $\gamma$  & $t=h^\vee$\\   
\hline    
$A_n$ &  $\mathfrak {sl}_{n+1}$     & $-2$ & 2 & $(n+1) $ & $n+1$\\
$B_n$ &   $\mathfrak {so}_{2n+1}$    & $-2$ & 4& $2n-3 $ & $2n-1$\\
$C_n$ & $ \mathfrak {sp}_{2n}$    & $-2$ & 1 & $n+2 $ & $n+1$\\
$D_n$ &   $\mathfrak {so}_{2n}$    & $-2$ & 4 & $2n-4$ & $2n-2$\\
$G_2$ &  $\mathfrak {g}_{2}  $    & $-2$ & $10/3 $& $8/3$ & $4$ \\
$F_4$ & $\mathfrak {f}_{4}  $    & $-2$ & $ 5$& $ 6$ & $9$\\
$E_6$ &  $\mathfrak {e}_{6}  $    & $-2$ & $ 6$& $ 8$ & $12$\\
$E_7$ & $\mathfrak {e}_{7}  $    & $-2$ & $ 8$& $ 12$ & $18$ \\
$E_8$ & $\mathfrak {e}_{8}  $    & $-2$ & $ 12$& $20$ & $30$\\
\hline  
\end{tabular}\label{tab:V1}
\end{table}

\begin{table}[h] 
\caption{Vogel's parameters for simple Lie algebras: lines}     
\begin{tabular}{|r|r|r|r|r|r|} 
\hline Algebra/Parameters & $\alpha$ &$\beta$  &$\gamma$  & $t$ & Line \\ 
\hline  $\mathfrak {sl}_{N}$  & -2 & 2 & $N$ & $N$ & $\alpha+\beta=0$ \\ 
\hline $\mathfrak {so}_{N}$ & -2  & 4 & $N-4$ & $N-2$ & $ 2\alpha+\beta=0$ \\ 
\hline  $ \mathfrak {sp}_{N}$ & -2  & 1 & $N/2+2$ & $N/2+1$ & $ \alpha +2\beta=0$ \\ 
\hline $Exc(n)$ & $-2$ & $2n+4$  & $n+4$ & $3n+6$ & $\gamma=2(\alpha+\beta)$\\ 
\hline 
\end{tabular} \caption*{For exceptional 
line $n=-2/3,0,1,2,4,8$ for $\mathfrak {g}_{2}, \mathfrak {so}_{8}, \mathfrak{f}_{4}, \mathfrak{e}_{6}, \mathfrak {e}_{7},\mathfrak {e}_{8} $, 
  respectively.} \label{tab:V2}
\end{table} 
Parameter $\alpha$ in these tables is chosen to be equal to $-2$. This always can be done due to the scaling invariance. 

For simple Lie algebras universal parameters have the following interpretation. Let's denote the eigenvalue of Casimir operator on adjoint representation as $2t$. Vogel \cite{V} shows that symmetric square of adjoint for all simple Lie algebras has the following  decomposition:  

\begin{eqnarray} \label{s2}
S^2 \mathfrak {g}=1+Y_2(\alpha)+Y_2(\beta)+Y_2(\gamma),\\
4t-2\alpha,4t-2\beta,4t-2\gamma
\end{eqnarray}
where second row contains  values of the same Casimir operator on representations $Y_2(\alpha),Y_2(\beta),Y_2(\gamma)$, respectively. This is actually definition of parameters $\alpha, \beta, \gamma$.  One can show that 
\begin{eqnarray}
\alpha+\beta+\gamma=t
\end{eqnarray}
and these are the same parameters introduced above, specialized for simple Lie algebras. Some subtlety in definition of irreducibility is in that one consider irreducibility w.r.t. the semidirect product of algebras on outer automorphisms of their Dynkin diagram, see \cite{Del,DM,Cohen}.  So, rescaling of parameters corresponds to rescaling of invariant scalar product in algebra, and permutation of parameters is equivalent to permutation of representations $Y_2(.)$. Choice $\alpha=-2$ corresponds to normalization of scalar product  such that the square of long roots is 2, and parameter $t$ is equal to dual Coxeter number $h^{\vee}$.

Unfortunately, abovementioned assumptions appear to be partially  wrong, particularly algebra $\Lambda$ appears to be not an integer domain \cite{V0}, so these simple modules actually are not correctly defined. Possible options are as following. 

First, the structure of algebra $\Lambda$ actually is not completely known. Conjecture of \cite{V} is that $\Lambda$ is generated by special diagrams $t$ and $x_n$, and in that case structure is known almost exactly. However, this conjecture is checked only in first few orders in a natural grading according to the number of vertexes, and it is perfectly possible that there are other generators in higher orders. If so, it may be that new generators appear at the same order when zero divisors relation appears, and some redefinition of product operation can lead to new ring which is already an integral domain. Unfortunately, these speculations are inaccessible for computer check, yet.

Second, although it is not possible to define simple modules with ring which is not integral domain, one can still have universal formulae. In that case they will have an arbitrariness of adding polynomials or other functions of universal parameters which are zero for all simple Lie algebras, and even on the entire lines of classical algebras, exceptional line and line $t=0$, which corresponds to the superalgebra $D_{2,1,\lambda}$ \cite{V0}.

There is a number of universal formulae purely in the theory of simple Lie algebras: universal formulae for dimensions of (series of) representations, appearing in decomposition of powers of adjoint representation \cite{LM1},  eigenvalues of higher Casimir operators \cite{V0,MSV}, volume of simple Lie groups \cite{M13,KhM},  quantum dimension of adjoint \cite{W3,MV}, etc. 

As mentioned, in the present paper we derive universal quantum dimensions for series of representations, namely for those for which Landsberg-Manivel \cite{LM1} obtain  universal dimension formulae. Particular cases of that formulae (for the cube of adjoint representation) were used in \cite{MMM} for derivation of universal invariant polynomials for two- and three-strands torus knots. As a check of these complicated formulae we prove numerically the particular case of Deligne's hypothesis on universal quantum dimensions, formulated below in Section \ref{s:4}.

Main application, discussed below in Section  \ref{s:1in},  is the universality of one-instanton  partition function of supersymmetric Yang-Mills theory.

\section{Quantum dimensions of Cartan powers of adjoint} \label{s:2}

As mentioned, we shall use  the theory of simple Lie algebras for calculation of universal quantum dimensions. Character of irreducible representation, restricted to the Weyl line (see definition below) is called below the universal character   of that representation. The Weyl formula for these characters gives (see e.g. \cite{DiF}, 13.170):

\begin{equation} \label{W}
\chi_{\lambda}(x\rho)= \prod_{\mu >0} \frac{\sinh(\frac{x}{2}(\mu,\lambda+\rho))}{\sinh(\frac{x}{2}(\mu,\rho))}
\end{equation}
Here $\rho$ is Weyl vector in root space, i.e. the half sum of positive roots; $\lambda$ is highest weight of irreducible representation; in l.h.s. we have the character of that representation, restricted to the Weyl line $x\rho$;  product is over positive roots $\mu$; (,) denote invariant (Cartan-Killing) scalar product in an arbitrary normalization.

The same object is often called quantum dimension and denoted

\begin{equation}
D_{\lambda} \equiv \chi_{\lambda}(x\rho)
\end{equation}

Note also that as an index of $D_{\lambda}$ and $\chi_{\lambda}$ we shall use also other notations, if they exist, for an irreducible representation, instead of its highest weight.

From (\ref{W}) we see that universal character is invariant w.r.t. the simultaneous rescaling of scalar product and $x$ variable: $x \rightarrow z x, (,) \rightarrow (,)/z$. Equivalently, universal formulae for quantum dimensions will be invariant w.r.t. the rescaling  $x \rightarrow z x, \alpha \rightarrow \alpha/z, \beta \rightarrow \beta/z, \gamma \rightarrow \gamma/z$ with an arbitrary $z$. Below we mainly use the normalization when square of long roots are equal 2 (this corresponds to $\alpha=-2$, \cite{LM1}). One can recover  parameter $\alpha$ by replacing $\beta \rightarrow -2\beta/\alpha, \gamma \rightarrow -2\gamma/\alpha, x\rightarrow -x\alpha/2$. 

Obviously, non-trivial contribution in (\ref{W}) give only roots $\mu$ with non-zero scalar product with $\lambda$. First we discuss  the case $\lambda=\theta$, where $\theta$ is the highest root of adjoint representation. We use the analysis of \cite{LM1}. With normalization of scalar product such that square of long roots is 2, roots with nonzero scalar product with highest weight $\theta$ are $\theta$ itself, with square 2, and a number of roots with scalar product with $\theta$ equal to 1. For these roots we need for (\ref{W}) their scalar products with $\rho$. These roots can be organized into three sequences with unit spacing between scalar products of consecutive elements  with $\rho$. These sequences are the following (we assume $\gamma \geq \beta \geq 0$ and the rank of algebra at least three): first sequence of length $t-2$, starting at root $\mu$ with $(\rho,\mu)=1$ and ending at some root $\nu$ with $(\rho,\nu)=t-2$, second sequence, of length $\gamma -2$,  starting at some root  $\mu$ with $(\rho,\mu)=\beta/2$ and ending at some root $\nu$ with $(\rho,\nu)= \gamma + \beta/2 -3$, and finally third sequence, of length $\beta -2$,  starting at some root  $\mu$ with $(\rho,\mu)=\gamma/2$ and ending at some root $\nu$ with $(\rho,\nu)= \beta + \gamma/2 -3$. 

We illustrate this picture for $E_7$ algebra in Table \ref{tab:E7-black}.

\begin{table}[h]
\caption{$(\rho,\mu)$ for all roots with $(\theta,\mu)=1$ for $E_7$}
	\centering
		\begin{tabular}{l|r|r|r|r|r|r|r|r|r|r|r|r|r|r|r|r}
		$(\rho,\mu)=$&1&2&3&4&5&6&7&8&9&10&11&12&13&14&15&16\\
				I &1&1&1&1&1&1&1&1&1&1&1&1&1&1&1&1\\
				II & & & &1&1&1&1&1&1&1&1&1&1& & & \\
				III & & & & & &1&1&1&1&1&1& & & & & \\ 
		
		\end{tabular}
		\label{tab:E7-black}
		\end{table}

From this classification of roots one can easily calculate the quantum dimension of adjoint representation by Weyl formula (\ref{W}). One factor is that from highest weight of adjoint, $\theta$:

\begin{eqnarray}\label{Ff}
F(x,\alpha,\beta,\gamma)=\frac{\text{Sinh}\left[\frac{(2\gamma+2\beta+\alpha)x}{4}\right]}{\text{Sinh}\left[\frac{(2\gamma+2\beta+3\alpha)x}{4}\right]}
\end{eqnarray}

 Another contribution is from the roots with scalar product with $\theta$, equal to 1. One observe that consecutive numerators and denominators in each string in Weyl formula cancel, and remaining three border terms are: 
 
\begin{eqnarray} \label{Bf}
B(x,\alpha,\beta,\gamma)= -\frac{\text{Sinh}\left[\frac{(\gamma+2\beta+2\alpha)x}{4}\right]}{\text{Sinh}\left[\frac{\gamma x}{4}\right]}\frac{\text{Sinh}\left[\frac{(2\gamma+\beta+2\alpha)x}{4}\right]}{\text{Sinh}\left[\frac{\beta x}{4}\right]}\frac{\text{Sinh}\left[\frac{(2\gamma+2\beta+3\alpha)x}{4}\right]}{\text{Sinh}\left[\frac{\alpha x}{4}\right]}
\end{eqnarray}
Altogether we get for quantum dimension:

\begin{eqnarray} \nonumber \label{cad}
f(x) &=& -\frac{\text{Sinh}\left[\frac{(\gamma+2\beta+2\alpha)x}{4}\right]}{\text{Sinh}\left[\frac{\gamma x}{4}\right]}\frac{\text{Sinh}\left[\frac{(2\gamma+\beta+2\alpha)x}{4}\right]}{\text{Sinh}\left[\frac{\beta x}{4}\right]}\frac{\text{Sinh}\left[\frac{(2\gamma+2\beta+\alpha)x}{4}\right]}{\text{Sinh}\left[\frac{\alpha x}{4}\right]}
\end{eqnarray}

Note that in the final answer Vogel's parameters enter symmetrically, as it should be for adjoint representation.

Similarly, quantum dimension of n-th Cartan power of adjoint is the product of contribution in Weyl formula of highest weight of adjoint representation ($F(x,n,\alpha,\beta,\gamma)$ below) and contribution of other roots ($B(x,n,\alpha,\beta,\gamma)$ below):

\begin{eqnarray} \label{nad}
& \chi_{n\theta}(x\rho)= F(x,n,\alpha,\beta,\gamma)B(x,n,\alpha,\beta,\gamma) \\ 
& F(x,n,\alpha,\beta,\gamma)=\frac{\text{Sinh}\left[\frac{(2\gamma+2\beta-(-3+2n)\alpha)x}{4}\right]}{\text{Sinh}\left[\frac{(2\gamma+2\beta+3\alpha)x}{4}\right]}  \\
& B(x,n,\alpha,\beta,\gamma)= \\ \nonumber
& \prod _{i=1}^n \frac{\text{Sinh}\left[\frac{(\gamma+2\beta-(-3+i)\alpha)x}{4}\right]}{\text{Sinh}\left[\frac{(\gamma-(-1+i)\alpha)x}{4}\right]}\frac{\text{Sinh}\left[\frac{(2\gamma+\beta-(-3+i)\alpha)x}{4}\right]}{\text{Sinh}\left[\frac{(\beta-(-1+i)\alpha)x}{4}\right]}\frac{\text{Sinh}\left[\frac{(2\gamma+2\beta-(-4+i)\alpha)x}{4}\right]}{\text{Sinh}\left[\frac{-i \alpha x}{4}\right]}
 \end{eqnarray}

For n=1 these functions coincide with those for adjoint representation (\ref{Ff}), (\ref{Bf}), (\ref{cad}):

\begin{eqnarray}
F(x,1,\alpha,\beta,\gamma)&=& F(x,\alpha,\beta,\gamma) \\
B(x,1,\alpha,\beta,\gamma)&=& B(x,\alpha,\beta,\gamma) \\
\chi_{n\theta}(x\rho)|_{n=1}&=& f(x)
\end{eqnarray}

For $n=2, 3$ in the limit $x\rightarrow 0$ we get dimensions of representations $Y_2(\alpha), Y_3(\alpha)$  coinciding with expressions of \cite{V}. For general $n$ in the same limit we obtain dimension formula for $Y_n(\alpha)$ from \cite{LM1}.

Although above classification of roots is valid not for all algebras, final result for characters is valid for all simple Lie algebras, as we checked directly. E.g. for $G_2$ we use Weyl formula (\ref{W}) with following $G_2$ data: simple roots $\alpha_1$ and $\alpha_2$, with length squares 2/3 and 2, respectively, and scalar product $(\alpha_1,\alpha_2)=-1$; remaining positive roots are $\alpha_1+\alpha_2$, $2\alpha_1+\alpha_2$, $3\alpha_1+\alpha_2$ and  $3\alpha_1+2\alpha_2$, last one being highest root. Substituting this in (\ref{W}) with $\lambda=n (3\alpha_1+2\alpha_2)$ we get 

\begin{eqnarray}
& \frac{\text{Sinh}\left[\frac{x}{2}(n+1)\right]}{\text{Sinh}\left[\frac{x}{2}\right]}\frac{\text{Sinh}\left[\frac{x}{2}\left(n+\frac{4}{3}\right)\right]}{\text{Sinh}\left[\frac{2x}{3}\right]}\frac{\text{Sinh}\left[\frac{x}{2}\left(n+\frac{5}{3}\right)\right]}{\text{Sinh}\left[\frac{5x}{6}\right]}\frac{\text{Sinh}\left[\frac{x}{2}(n+2)\right]}{\text{Sinh}[x]}\frac{\text{Sinh}\left[\frac{x}{2}(2n+3)\right]}{\text{Sinh}\left[\frac{3x}{2}\right]}
\end{eqnarray}
which coincides with general formula (\ref{nad}) with  universal parameters specialized for $\mathfrak {g}_{2}$.

\section{Quantum dimensions of Cartan powers of adjoint and $Y_2(\beta)$} \label{s:3}

Next consider the Cartan product of $k$ copies of adjoint representation and $l$ copies of  representation  $Y_2(\beta)$ . We again follow the analysis of \cite{LM1}.

Consider highest root $\theta$ of simple Lie algebra, corresponding principal $sl_2$, its centralizer $\mathfrak{h}$ and $\mathfrak{h}$'s highest root $\sigma$, assuming $\mathfrak{h}$ is simple Lie algebra. See Table (\cite{LM1}, p.385) for a list of corresponding data. According to \cite{LM1,M1}, if $\mathfrak{h}$ is simple, its  Vogel's parameters $\alpha', \beta', \gamma'$  are $\alpha'=\alpha, \beta'=\gamma -\beta, \gamma'=\beta$, provided initial parameters are ordered as $\alpha < 0, \gamma \geq \beta \geq 0$.  Since $\theta,\sigma $ are orthogonal, their corresponding Cartan elements $H_{\theta}, H_{\sigma}$ create a double grading 

\begin{eqnarray}
g_{ij}=\{X\in g, [H_{\theta},X]=iX, [H_{\tilde{\alpha}},X]=jX \}
\end{eqnarray}  

Non-zero are $g_{00},g_{0\pm 1}, g_{\pm 10}, g_{\pm 1 \pm 1} , g_{\pm 2 0}, g_{0\pm 2}$ spaces. We assume minimal normalization of invariant scalar product, i.e. square of long roots =2, which corresponds to Vogel's parameter $\alpha=-2$.

The highest weight of $Y_2(\beta)$ is $\theta + \sigma$ \cite{LM1}, where $\sigma$ is the highest weight of subgroup $\mathfrak{h}$, so the highest weight of Cartan product of $k$ copies of adjoint and $l$ copies of $Y_2(\beta)$ is $(k+l)\theta + l \sigma$.

For Weyl formula we need  $g_{01},g_{0 2},g_{2 0},  g_{10}, g_{1 \pm 1} $.  Last two spaces, or equivalently  roots with scalar product with highest root of $\mathfrak{g}$ equal to 1, are already used in previous section for calculation of quantum dimension of powers of adjoint representation and are presented in Table \ref{tab:E7-black}  for $E_7$ as an example. Universal parameters of $E_7$ are $(\alpha,\beta,\gamma)=(-2,8,12)$.

\begin{table}[h]
\caption{$(\rho,\mu)$ for roots $\mu$ of $g_{1\pm 1}$ ((r)ed and (g)reen) and $g_{10}$ (black) for $E_7$}
	\centering
		\begin{tabular}{l|r|r|r|r|r|r|r|r|r|r|r|r|r|r|r|r}
		$(\rho,\mu)=$&1&2&3&4&5&6&7&8&9&10&11&12&13&14&15&16\\
				I & r & r & r & r & r & r & r & 1 & 1 & g & g & g & g & g & g & g \\
				II & & & & r & 1 & 1 & 1 & 1 & 1 & 1 & 1 & 1 & g & & & \\
				III & & & & & & 1 & 1 & 1 & 1 & 1 & 1& & & & & \\ 
		
		\end{tabular}
		\label{tab:E7Roots}
		\caption*{Roots $\mu$ with $(\theta,\mu)=1$. Black roots have $(\sigma,\mu)=0$, (r)ed roots have $(\sigma,\mu)=-1$, (g)reen roots  have $(\sigma,\mu)=1$. All entries have multiplicity 1. }
\end{table}

Now we need refinement of space of roots with unit scalar product with $\theta$ according to their scalar product with $\sigma$. These are shown in the  Table \ref{tab:E7Roots} for $E_7$, and described for general case below. Red and green roots belong to $g_{1,-1}$  and $g_{1,1}$, respectively. 

The number of both is $\beta$, with $(\beta-1)$ at the beginning (end) of sequence I, plus one additional root at the beginning (end)  of sequence II. 

Now we have all necessary data for Weyl formula. First we  write down  contributions  from  different parts of three sequences of roots.  Red roots  give $C_2(x,k,\alpha,\beta,\gamma)$, green ones $C_1(x,k+2l,\alpha,\beta,\gamma)$, remaining colored roots give $F_{1-1},F_{11}$ respectively. Three sequences of remaining black roots give $A(x,k+l,\alpha,\beta,\gamma)$. Next, one-dimensional spaces  $g_{20}, g_{02}$ give rise to $F_{20},F_{02}$, respectively. Finally, contribution of $g_{01}$ we denote by $\tilde{B}(x,l,\alpha,\beta,\gamma)$. Contribution of  $g_{01}$ and $g_{02}$ together give the quantum dimension of $l$-th power of adjoint of $\mathfrak h$. Since we know universal parameters of $\mathfrak h$, we can use calculation of quantum dimension of power of adjoint representation in previous section. Particularly, $\tilde{B}(x,l,\alpha,\beta,\gamma) = B(x,l,\alpha,\beta,\gamma-\beta)$.

Altogether,  the product of all contributions gives desired answer $\chi_{(k+l)\theta+l\sigma}$ which we denote  $Z(x,k,l,\alpha,\beta,\gamma)$

\begin{eqnarray}
\chi_{(k+l)\theta+l\sigma}(x\rho)\equiv Z(x,k,l,\alpha,\beta,\gamma)
\end{eqnarray} 

Contribution  of black roots in Table \ref{tab:E7Roots} into character (for an arbitrary simple group) is given by  $A(x,n,\alpha,\beta,\gamma)$ below at $n=k+l$:
\begin{eqnarray} \nonumber
& 
A(x,n,\alpha,\beta,\gamma)=\prod _{i=1}^n \left(\frac{\text{Sinh}\left[\frac{(2\gamma+3\alpha-i \alpha)x}{4}\right]}{\text{Sinh}\left[\frac{(2\beta+\alpha-i \alpha)x}{4}\right]}\frac{\text{Sinh}\left[\frac{(2\gamma+\beta+4\alpha-i \alpha)x}{4}\right]}{\text{Sinh}\left[\frac{(\beta-i \alpha)x}{4}\right]}\frac{\text{Sinh}\left[\frac{(2\beta+\gamma+3\alpha-i \alpha)x}{4}\right]}{\text{Sinh}\left[\frac{(\gamma+\alpha-i \alpha)x}{4}\right]}\right)
\end{eqnarray}

Contribution of red roots is given by function below at $n=k$:
\begin{eqnarray}
C_2(x,n,\alpha,\beta,\gamma)=\prod _{i=1}^n \frac{\text{Sinh}\left[\frac{x}{4} (-\alpha-2 \beta+\alpha i)\right]}{\text{Sinh}\left[\frac{\alpha i x}{4}\right]}
\end{eqnarray}

Contribution of green roots is given by function below at $n=k+2l$:
\begin{eqnarray}
C_1(x,n,\alpha,\beta,\gamma)=\prod _{i=1}^n \frac{\text{Sinh}\left[\frac{(2\beta+2\gamma+4\alpha-i \alpha)x}{4}\right]}{\text{Sinh}\left[\frac{(2\gamma+3\alpha-i \alpha)x}{4}\right]}
\end{eqnarray}

Contribution of 1-dimensional spaces $g_{20},g_{02}$,  and remaining colored roots from $g_{11}, g_{1-1}$:
\begin{eqnarray} \nonumber
& 
\text{F}(x,k,l,\alpha,\beta,\gamma) \text{=} F_{20}F_{02}F_{11}F_{1,-1} = \\ \nonumber 
& = \frac{\text{Sinh}\left[\frac{(2\beta+2\gamma-(-3+2k+2l)\alpha)x}{4}\right]}{\text{Sinh}\left[\frac{(2\beta+2\gamma+3\alpha)x}{4}\right]}\frac{\text{Sinh}\left[\frac{(2\gamma-(-3+2l)\alpha)x}{4}\right]}{\text{Sinh}\left[\frac{(2\gamma+3\alpha)x}{4}\right]}\frac{\text{Sinh}\left[\frac{(\beta+2\gamma-(-3+k+2l)\alpha)x}{4}\right]}{\text{Sinh}\left[\frac{(\beta+2\gamma+3\alpha)x}{4}\right]}\frac{\text{Sinh}\left[\frac{(\beta-k \alpha)x}{4}\right]}{\text{Sinh}\left[\frac{\beta x}{4}\right]}
\end{eqnarray}

Contribution of $g_{01}$:
\begin{eqnarray} \nonumber
&\tilde{B}(x,l,\alpha,\beta,\gamma) = B(x,l,\alpha,\beta,\gamma-\beta)= \\
& \prod _{i=1}^n \frac{\text{Sinh}\left[\frac{(2\gamma-\beta-(-3+i)\alpha)x}{4}\right]}{\text{Sinh}\left[\frac{(\beta-(-1+i)\alpha)x}{4}\right]}\frac{\text{Sinh}\left[\frac{(\gamma+\beta-(-3+i)\alpha)x}{4}\right]}{\text{Sinh}\left[\frac{(\gamma-\beta-(-1+i)\alpha)x}{4}\right]}\frac{\text{Sinh}\left[\frac{(2\gamma-(-4+i)\alpha)x}{4}\right]}{\text{Sinh}\left[\frac{-i \alpha x}{4}\right]}
\end{eqnarray}

Finally, quantum dimension of Cartan product $g^k Y_2(\beta)^l$ is:
\begin{eqnarray}\nonumber
& Z(x,k,l,\alpha,\beta,\gamma)= \\
& F(x,k,l,\alpha,\beta,\gamma)A(x,k+l,\alpha,\beta,\gamma) \times \\ 
&\tilde{B}(x,l,\alpha,\beta,\gamma)C_1(x,k+2l,\alpha,\beta,\gamma)C_2(x,k,\alpha,\beta,\gamma)
\end{eqnarray}
In the limit $x \rightarrow 0$ one obtains universal formulae for dimensions of corresponding representation. They coincide with formulae of \cite{LM1}.

Again, we check this expression for algebras, for which derivation is not valid, directly.  However, there is an interesting subtlety. Let's take again example of rank-two algebra $G_2$. One can calculate the quantum dimension for Cartan  product of $k$ adjoint representations and $p$ representations $Y_2(\beta)$ using data for $G_2$ given above:

\begin{eqnarray}
& f(x,k,p)=\frac{\text{Sinh}\left[\frac{x}{2}\left(\frac{2p+1}{3}\right)\right]}{\text{Sinh}\left[\frac{x}{2}\left(\frac{1}{3}\right)\right]}\frac{\text{Sinh}\left[\frac{x}{2}(k+1)\right]}{\text{Sinh}\left[\frac{x}{2}(1)\right]}\frac{\text{Sinh}\left[\frac{x}{2}\left(k+1+\frac{2p+1}{3}\right)\right]}{\text{Sinh}\left[\frac{x}{2}\left(1+\frac{1}{3}\right)\right]} \times \\ \nonumber
& \frac{\text{Sinh}\left[\frac{x}{2}\left(k+1+\frac{2(2p+1)}{3}\right)\right]}{\text{Sinh}\left[\frac{x}{2}\left(1+\frac{2}{3}\right)\right]}\frac{\text{Sinh}\left[\frac{x}{2}\left(k+1+\frac{3(2p+1)}{3}\right)\right]}{\text{Sinh}\left[\frac{x}{2}\left(1+\frac{3}{3}\right)\right]}\frac{\text{Sinh}\left[\frac{x}{2}\left(2(k+1)+\frac{3(2p+1)}{3}\right)\right]}{\text{Sinh}\left[\frac{x}{2}\left(2+\frac{3}{3}\right)\right]}
\end{eqnarray}

One can easily check coincidence for $p=1$:

\begin{eqnarray}
Z(x,k,1,-2,10/3,8/3)=f(x,k,1)
\end{eqnarray}

However, we encounter zero for  $p>1$:

\begin{eqnarray} \label{z=0}
Z(x,k,p,-2,10/3,8/3)=0 \\
\text{for} \,\, p>1
\end{eqnarray}

It is not clear why is it so, but it is in agreement with direct calculations at low orders for exceptional line in \cite{Cohen}. Namely, it was found there that for $k=0, p=2$ the corresponding universal representation $J$ (universal in the sense that there exist universal - on the exceptional line - expression for its dimension, see below)  is equal to Cartan square of  $Y_2(\beta)$ (which is denoted there $Y^*$) for all exceptional algebras except $\mathfrak {g}_{2}$, for which it is equal to zero. Our representation with $k=0,p=2$ coincides with $J$, as can be seen either from decomposition formulae from \cite{Cohen}, or one can obtain dimension formula  (restricting  $Z(x,0,2,\alpha,\beta,\gamma)$ to exceptional line and taking limit $x \rightarrow 0$), and then compare the result with dimension formula for $J$. Parameterizing exceptional line by  $\alpha=\lambda,\beta=1-\lambda,\gamma=2$ and taking the limit we get dimension formula

\begin{eqnarray}
& \lim_{x\rightarrow 0}Z(x,0,2,\lambda,1-\lambda,2)= \\ \nonumber
& \frac{81 (-6+\lambda ) (-4+\lambda ) (-3+\lambda ) (2+\lambda ) (3+\lambda ) (5+\lambda ) (-5+2 \lambda ) (3+2 \lambda )}{(-1+\lambda )^2 \lambda ^2 (-1+2 \lambda )^2 (-2+3 \lambda ) (-1+3 \lambda )}
\end{eqnarray}
which coincides exactly with dimension formula for representation $J$ from \cite{Cohen}. (Similar coincidence one can observe with other formulae from \cite{Cohen}, see table in Appendix.) In this parametrization $\mathfrak {g}_{2}$ corresponds to $\lambda=-2/3$, and above expression is zero at that point.

So, phenomenon (\ref{z=0}) is an extension of this fact to the whole $k\geq 0, p\geq 2$ region. This and other specific features of universal formulae for quantum dimensions deserve further study.

\section{Deligne's hypothesis on quantum dimensions} \label{s:4}

Introduction of parameters $\alpha, \beta, \gamma$ can be considered as some deformation of the scalar objects in the theory of (simple) Lie algebras and groups. The reasonable question is what properties of initial objects are maintained.  Deligne \cite{D13} suggested that usual relation between characters of representations (i.e. the product of characters of representations is equal to the sum of characters of representations in decomposition of product of representations) is satisfied by universal quantum dimensions at entire Vogel's plane. This hypothesis would follow from Vogel's Universal Lie algebra \cite{V} provided, first, one defines the quantum dimensions in the framework of Universal Lie algebra and, second, there exist a homomorphism of Universal Lie algebra into complex numbers with an arbitrary parameters $\alpha, \beta, \gamma$. Both conditions are not established yet.

Deligne carried on  complete check of the version of his hypothesis restricted to  $\mathfrak{sl}_{N}$ line \cite{D13}. On the entire Vogel's plane one can check it for the square  of adjoint representation. For example, the symmetric square of adjoint representation of simple Lie algebra has a universal decomposition \cite{V0}, given above:

\begin{eqnarray} \label{s22}
S^2 \mathfrak {g}=1+Y_2(\alpha)+Y_2(\beta)+Y_2(\gamma)
\end{eqnarray}

Quantum dimension of  representation $ Y_2(\alpha)$ is given in Section \ref{s:2}:

\begin{eqnarray} \label{Y2c}
& \chi_{Y_2(\alpha)}(x\rho) =  \\ \nonumber
&  -\frac{\text{sinh}[\frac{x t}{2}]\text{sinh}[\frac{x(\beta-2t)}{4}]\text{sinh}[\frac{x(\gamma-2t)}{4}]\text{sinh}[\frac{x(\beta+t)}{4}]\text{sinh}[\frac{x(\gamma+t)}{4}]\text{sinh}[\frac{x(3\alpha-2t)}{4}]}{\text{sinh}[\frac{x\alpha}{4}]\text{sinh}[\frac{x\alpha}{2}]\text{sinh}[\frac{x\beta}{4}]\text{sinh}[\frac{x\gamma}{4}]\text{sinh}[\frac{x(\alpha-\beta)}{4}]\text{sinh}[\frac{x(\alpha-\gamma)}{4}]}
\end{eqnarray}
and permutations of this for $Y_2(\beta),Y_2(\gamma)$. 

Deligne's hypothesis imply 

\begin{eqnarray} \label{2sym}
 \label{decompchar}
 \chi_{S^2\mathfrak{g}}(x\rho)=\chi_{Y_2(\alpha)}(x\rho)+\chi_{Y_2(\beta)}(x\rho)+\chi_{Y_2(\gamma)}(x\rho)+ 1
 \end{eqnarray}
at an arbitrary universal parameters. Here quantum dimension of symmetric square of adjoint can be expressed in terms of quantum dimension of adjoint, i.e. function $f(x)$ (\ref{cad}):

\begin{eqnarray}\label{2sym2}
\chi_{S^2\mathfrak{g}}(x\rho)=\frac{1}{2} \left( f^2(x)+f(2x) \right)
\end{eqnarray}

Then relation (\ref{2sym}) can be checked directly. 
Similarly, decomposition of antisymmetric square of adjoint:

\begin{eqnarray}
 \label{decompchar2}
 \wedge^2(\mathfrak{g})= \mathfrak{g}+X_2
 \end{eqnarray}

implies relation between quantum dimensions:

\begin{eqnarray}
 \label{decompchar2-2}
\frac{1}{2} \left( f^2(x)-f(2x) \right)= f(x)+ D_{X_2}
 \end{eqnarray}

where quantum dimension of representation $X_2$ is \cite{D13}

\begin{eqnarray} 
&
D_{X_2}=\frac{\text{Sinh}\left[\frac{(2t-\alpha)x}{4}\right]\text{Sinh}\left[\frac{(2t-\beta)x}{4}\right]\text{Sinh}\left[\frac{(2t-\gamma)x}{4}\right]}{\text{Sinh}\left[\frac{\alpha x}{4}\right]\text{Sinh}\left[\frac{\beta x}{4}\right]\text{Sinh}\left[\frac{\gamma x}{4}\right]}  \times \\ \nonumber
&
\frac{\text{Sinh}\left[\frac{(t+\alpha)x}{4}\right]\text{Sinh}\left[\frac{(t+\beta)x}{4}\right]\text{Sinh}\left[\frac{(t+\gamma)x}{4}\right]}{\text{Sinh}\left[\frac{\alpha x}{2}\right]\text{Sinh}\left[\frac{\beta x}{2}\right]\text{Sinh}\left[\frac{\gamma x}{2}\right]}   
\frac{\text{Sinh}\left[\frac{(t-\alpha)x}{2}\right]\text{Sinh}\left[\frac{(t-\beta)x}{2}\right]\text{Sinh}\left[\frac{(t-\gamma)x}{2}\right]}{\text{Sinh}\left[\frac{(t-\alpha)x}{4}\right]\text{Sinh}\left[\frac{(t-\beta)x}{4}\right]\text{Sinh}\left[\frac{(t-\gamma)x}{4}\right]}
\end{eqnarray}

Now we would like to check the hypothesis for symmetric cube of adjoint representation. 

Decomposition of symmetric cube of adjoint is, according to \cite{V}:

\begin{eqnarray} \label{s3}
& S^3\mathfrak{g}= 2 \mathfrak{g} +X_2+Y_3(\alpha)+Y_3(\beta)+Y_3(\gamma)+ \\ \nonumber
&\mathfrak{g}Y_2(\beta)(\alpha,\beta,\gamma)+\mathfrak{g}Y_2(\beta)(\alpha,\gamma,\beta)+\mathfrak{g}Y_2(\beta)(\beta,\gamma,\alpha)
\end{eqnarray}
where $\mathfrak{g}Y_2(\beta)(\alpha,\beta,\gamma)$ is the Cartan  product of adjoint and $Y_2(\beta)$ representations, where we explicitly write down its dependence on universal parameters, and last two representations are obtained from $\mathfrak{g}Y_2(\beta)(\alpha,\beta,\gamma)$ by permutations of parameters. Note that   $Z(x,1,1,\alpha,\beta,\gamma)$ is symmetric w.r.t. the switch of $\alpha$ and $\beta$, so there is no other representations. In Appendix we give a specific examples of this decomposition. 

Equation (\ref{s3}) leads to the probable relation between quantum dimensions:

\begin{eqnarray} \label{40} \nonumber
& & \frac{1}{6}\left(f(x,\alpha,\beta,\gamma)^3+3f(2x,\alpha,\beta,\gamma)f(x,\alpha,\beta,\gamma)+2f(3x,\alpha,\beta,\gamma)\right)= \\
& & Z(x,3,0,\alpha,\beta,\gamma)+Z(x,3,0,\beta,\alpha,\gamma)+Z(x,3,0,\gamma,\beta,\alpha)+ \\ \nonumber 
& & Z(x,1,1,\alpha,\beta,\gamma)+  Z(x,1,1,\alpha,\gamma,\beta)+Z(x,1,1,\beta,\gamma,\alpha)+\\ \nonumber 
& & +D_{X_2}(x,\alpha,\beta,\gamma)+2f(x,\alpha,\beta,\gamma)
\end{eqnarray}
which is obliged to be satisfied on the points of Vogel's table, but Deligne's hypothesis assumes it on entire Vogel's plane. 
This is checked with Mathematica  up to 17-th order in expansion over x, and numerically at a dozens thousands of random points in Vogel's plane and $x$ line. All this cannot be considered as a rigorous proof, but is convincing enough and can be developed into such. 

 Finally, let's note that universal quantum dimensions, entering in decomposition of symmetric cube of adjoint representation, i.e. those in r.h.s. of (\ref{40}), were used in \cite{MMM} in calculation of universal knot polynomials.

\section{One-instanton partition function} \label{s:1in}

Universal formulae of present paper can have a number of applications. We would like to mention one-instanton contribution into Nekrasov's partition function for an arbitrary group in pure N=2 4d superYang-Mills theory. 

It is calculated for an arbitrary gauge group $G$ in \cite{KMJT}, see also \cite{BHM}. It is essentially given (see Appendix B of \cite{KMJT}) by $\sigma \rightarrow 0$ limit of character  of representation

\begin{eqnarray} \label{1in}
\sum_{n=1}^{\infty} V(-n\theta)\otimes T^{\otimes n}
\end{eqnarray}
of group $G \otimes U(1)^2$. Here $V(-n\theta)$ is an irrep with highest weight $n\theta$, i.e. $n$-th Cartan power of adjoint, $T=T_1T_2$, $T_i=\exp(\sigma \epsilon_i), \, i=1,2$, $\epsilon_i$ are Nekrasov's parameters. 

Universalization of this expression is possible, if vacuum expectation value of scalar field (argument of character) is restricted to Weyl line $x\rho$. Then in formula (\ref{1in}) appears the universal character of $n$-th Cartan power of adjoint given by (\ref{nad}), and we get a universal expression for 1-instanton contribution:  

\begin{eqnarray} \label{1in2} 
& \sum_{n=1}^{\infty}e^{n\sigma(\epsilon_1+\epsilon_2)}\frac{\text{Sinh}\left[\frac{(2\gamma+2\beta-(-3+2n)\alpha)x}{4}\right]}{\text{Sinh}\left[\frac{(2\gamma+2\beta+3\alpha)x}{4}\right]} \times \\ \nonumber
&
\prod _{i=1}^n \frac{\text{Sinh}\left[\frac{(\gamma+2\beta-(-3+i)\alpha)x}{4}\right]}{\text{Sinh}\left[\frac{(\gamma-(-1+i)\alpha)x}{4}\right]}\frac{\text{Sinh}\left[\frac{(2\gamma+\beta-(-3+i)\alpha)x}{4}\right]}{\text{Sinh}\left[\frac{(\beta-(-1+i)\alpha)x}{4}\right]}\frac{\text{Sinh}\left[\frac{(2\gamma+2\beta-(-4+i)\alpha)x}{4}\right]}{\text{Sinh}\left[\frac{-i \alpha x}{4}\right]}
 \end{eqnarray}

This observation can be developed in the several directions: one can try to get closed form for the sum over $n$ of universal characters in (\ref{1in}), extend this to higher-instanton contributions and  to complete partition function of (supersymmetric) Yang-Mills theory, and to other ("universal") values of scalar fields. 

\section{Conclusion}

As mentioned, this paper is the next step in realization of program of representation of partition functions and observables in gauge theories (and other simple-Lie-algebras-based theories) in universal form. 
This program is successful for Chern-Simons theory on the 3d sphere \cite{MV,M13}, and has led to establishment of exact (i.e. non-perturbative) Chern-Simons/topological strings duality for SU(N) gauge groups \cite{KM}, similar results in refined cases \cite{KM}, gauge/string duality conjecture for exceptional groups \cite{M14}, universal knot polynomials \cite{MMM}, and other achievements.

In present paper we derive the universal expressions for universal characters (quantum dimensions) of some series of representation of simple Lie algebras, appearing in decomposition of tensor powers of adjoint representation. We partially  check these complicated formulae by numerical study of particular case of Deligne hypothesis on quantum dimensions. 

There is a number of applications of these formulae for abovementioned program. Our expressions for quantum dimensions themselves are quantum averages for unknot Wilson loops, in corresponding representations, in Chern-Simons theory on 3d sphere. Next, some of these formulae, namely those for irreps in the decomposition of the cube of adjoint representation are already used in \cite{MMM} for derivation of universal form of knot polynomials (i.e. Wilson averages) for some torus knots in the same Chern-Simons theory. In  Section 5 we use these formulae for derivation of universal form of one-instanton partition function on Weyl line in supersymmetric Yang-Mills theory.

We expect further applications of present formulae in gauge theories. Most important would be derivation of universal form of perturbative and non-perturbative parts of (supersymmetric) Yang-Mills theories  on different manifolds and with different matter fields (super)multiplets. With the same theories in mind it seems intriguing to discover universality in integrable models.

\section{Acknowledgments}
I'm indebted to P.Deligne for discussions and correspondence at autumn of 2013, when  part of this work has been done, and to R.Poghossian for discussions on instanton sums.

Work is partially supported by the Science Committee of the Ministry of Science 
and Education of the Republic of Armenia under contract  15T-1C233.

\appendix

\section{Examples of decompositions of symmetric cube of adjoint} \label{s:5}

In this Appendix we present three Tables for decompositions of symmetric cube of adjoint representation for specific algebras: $\mathfrak {sl}_{6} , \mathfrak {sl}_{N},  \mathfrak {so}_{12}, \mathfrak {so}_{N}, \mathfrak {f}_{4}$ and exceptional line. We see particularly appearance of "negative" representations (see data for $\mathfrak {f}_{4}$) which actually cancel the same representation, appearing in some other term in decomposition, so effectively they don't appear at all.

\begin{table}[ht]  
\caption{Decomposition of the symmetric cube of adjoint for $\mathfrak {sl}_{6}$ and $\mathfrak {sl}_{N}$} 
\begin{tabular}{|c|c|c|c|}
\hline irrep & dim ($\mathfrak {sl}_{6}$) &    Dynkin  & dim($\mathfrak {sl}_{N}$)  \\ 
\hline $2\mathfrak{g}$ & $2\times 35$ &  $2\times$ 10001 & $2 (-1 + N^2)$  \\ 
\hline $Y_3(\alpha)$ & 2695 &  30003 & $\frac{1}{36} (-1+N) N^2 (1+N)^2 (5+N)$  \\ 
\hline $Y_3(\beta)$ & 175  &   00200 & $\frac{1}{36} (-5+N) (-1+N)^2 N^2 (1+N)$   \\ 
\hline $Y_3(\gamma)$ & 1 &    & 1   \\ 
\hline $X_2$ & $2\times 280$  &   20010+01002 & $\frac{1}{2} (-2+N) (-1+N) (1+N) (2+N)$   \\ 
\hline $\mathfrak{g}Y_2(\beta)(\alpha,\beta,\gamma)$ & 3675 &   11011 & $\frac{1}{9} (-3+N) (3+N) \left(-1+N^2\right)^2 $   \\ 
\hline $\mathfrak{g}Y_2(\beta)(\alpha,\gamma,\beta)$ & 405 &   20002 & $\frac{1}{4} (-1+N) N^2 (3+N) $   \\ 
\hline $\mathfrak{g}Y_2(\beta)(\beta,\gamma,\alpha)$ & 189 &   01010  & $\frac{1}{4} (-3+N) N^2 (1+N)$  \\ 
\hline Sum of dims & 7770 &  & $\frac{1}{6} (-1+N) N^2 (1+N) \left(1+N^2\right)$ \\
\hline Dim of $S^3\mathfrak{g}$ & 7770 & & $\frac{1}{6} (-1+N) N^2 (1+N) \left(1+N^2\right)$  \\
\hline
\end{tabular} 
\label{tab:par1}
\end{table}

 \begin{table}[ht]  
 \caption{Decomposition of the symmetric cube of adjoint for $\mathfrak {f}_{4}$ and $Exc$ line} 
 \begin{tabular}{|c|c|c|c|}
 \hline irrep & dim  $(\mathfrak {f}_{4})$ &    Dynkin  & $Exc$ line: $\alpha=s,\beta=1-s, \gamma=2$  \\ 
 \hline $2\mathfrak{g}$  & 2$\times $52 &   2$\times $  1000 & $-\frac{4 (-6+s) (5+s)}{(-1+s) s} $   \\ 
 \hline $Y_3(\alpha)$ & 12376 &   3000 & $-\frac{10 (-6+s) (-5+s) (-4+s) (5+s) (-6+5 s)}{(-1+s)^2 s^3 (-1+2 s) (-1+3 s)} $   \\ 
 \hline $Y_3(\beta)$ & 273  &   0010 & $-\frac{10 (-6+s) (3+s) (4+s) (5+s) (1+5 s)}{(-1+s)^3 s^2 (-1+2 s) (-2+3 s)} $   \\ 
 \hline $Y_3(\gamma)$ & -52 &  &   $ \frac{2 (-6+s) (5+s)}{(-1+s) s}$    \\ 
 \hline $X_2$ & 1274  &   0100 & $\frac{5 (-6+s) (-4+s) (3+s) (5+s)}{(-1+s)^2 s^2} $   \\ 
 \hline $\mathfrak{g} Y_2(\beta)(\alpha,\beta,\gamma)$ & 10829  &   1002  & $-\frac{27 (-6+s) (-5+s) (-4+s) (3+s) (4+s) (5+s)}{(-1+s)^2 s^2 (-2+3 s) (-1+3 s)} $   \\ 
 \hline $ \mathfrak{g}Y_2(\beta), (\alpha,\gamma,\beta)$  & 0 &    & 0  \\ 
 \hline  $\mathfrak{g}Y_2(\beta)(\beta,\gamma,\alpha)$ & 0 &     & 0  \\ 
 \hline Sum of dims & 24804 &  &$\frac{20 (-6+s) (5+s) \left(-60-s+s^2\right)}{(-1+s)^3 s^3} $  \\
 \hline Dim of $S^3\mathfrak{g}$ & 24804 & & $\frac{20 (-6+s) (5+s) \left(-60-s+s^2\right)}{(-1+s)^3 s^3} $  \\
 \hline
 \end{tabular} 
 \label{tab:par2}
 \end{table}

 \begin{table}[ht]  
 \caption{Decomposition of the symmetric cube of adjoint for $\mathfrak {so}_{12}$ and $\mathfrak {so}_{N}$ } 
 \begin{tabular}{|c|c|c|c|}
 \hline irrep & dim $(\mathfrak {so}_{12})$ &    Dynkin  & dim $(\mathfrak {so}_{N})$  \\ 
 \hline $2\mathfrak{g}$ & 2$\times$66 &   2$\times$ 010000 & $ (-1+N) N $  \\ 
 \hline $Y_3(\alpha)$ & 23100 &   030000 &  $\frac{1}{144} (-3+N) (-2+N) (-1+N) (2+N) (3+N) (4+N)$   \\ 
 \hline $Y_3(\beta)$ & 924  &   000020+000002 & $\frac{1}{720} (-5+N) (-4+N) (-3+N) (-2+N) (-1+N) N $   \\ 
 \hline $Y_3(\gamma)$ & 0 &    & 0   \\ 
 \hline $X_2$ & 2079  &   101000 & $\frac{1}{8} (-3+N) (-1+N) N (2+N)$   \\ 
 \hline  $\mathfrak{g} Y_2(\beta)(\alpha,\beta,\gamma)$ & 21021  &   010100  & $\frac{1}{80} (-5+N) (-2+N) (-1+N) N (1+N) (2+N)$   \\ 
 \hline $ \mathfrak{g}Y_2(\beta), (\alpha,\gamma,\beta)$  & 2860 &   210000  & $ \frac{1}{8} (-2+N) (-1+N) (1+N) (4+N) $  \\ 
 \hline   $\mathfrak{g}Y_2(\beta)(\beta,\gamma,\alpha)$& 0 &    &    \\ 
 \hline Sum of dims & 50116 &  & $ \frac{1}{48} (-1+N) N \left(2-N+N^2\right) \left(4-N+N^2\right)$  \\
 \hline Dim of $S^3\mathfrak{g}$ & 50116 &  & $ \frac{1}{48} (-1+N) N \left(2-N+N^2\right) \left(4-N+N^2\right)$  \\
 \hline
 \end{tabular} 
 \label{tab:par3}
 \end{table}

\FloatBarrier

\end{document}